\documentclass{au}
\usepackage{amssymb, natbib, latexcad}

\theoremstyle{definition}
\newtheorem{definition}{Definition}[section]
\newtheorem{example}[definition]{Example}

\theoremstyle{plain}
\newtheorem{theorem}[definition]{Theorem}
\newtheorem{lemma}[definition]{Lemma}
\newtheorem{corollary}[definition]{Corollary}
\newtheorem{proposition}[definition]{Proposition}

\theoremstyle{remark}
\newtheorem{remark}[definition]{Remark}

\newcommand{\vm}[4]{
    \left(\begin{array}{cc}#1 & #2\\ #3 & #4 \end{array}\right)
}                                       
\newcommand{\paige}[1]{M^*(#1)}         
\newcommand{\cyclic}[1]{C_{#1}}         
\newcommand{\dpr}[2]{#1\cdot#2}         
\newcommand{\vpr}[2]{#1\times#2}        
\newcommand{\spn}[1]{\langle#1\rangle}  
\newcommand{\field}[1]{GF(#1)}          
\newcommand{\neutral}{e}                
\newcommand{\fullblank}{\_}             
\newcommand{\order}[1]{|#1|}            
\newcommand{\aut}[1]{\mathrm{Aut}(#1)}  
\newcommand{\chein}{M_{2n}(G,\,2)}      
\newcommand{\schein}{M(G)}              
\newcommand{\subalg}[1]{\mathrm{Sub}(#1)
}                                       
\newcommand{\join}{\vee}                
\newcommand{\meet}{\wedge}              
\newcommand{\nbd}[1]{\mathrm{Nbd}(#1)}  
\newcommand{\octo}[1]{O(#1)}            
\newcommand{\dswitch}{\partial}         
\newcommand{\ifix}{\xi}                 
\newcommand{\conj}[1]{\gamma_{#1}}      
\newcommand{\nucleus}[1]{N(#1)}         
\newcommand{\inv}[2]{[#1,\,#2]}         
\newcommand{\iinv}[3]{\inv{#1}{#2}_{#3}}
\newcommand{\tri}[3]{\{#1,\,#2\}_{#3}}  
\newcommand{\weight}[1]{\mathrm{w}(#1)} 
\newcommand{\seq}[1]{\stackrel{#1}{=}}  
\newcommand{\orbitname}{O}              
\newcommand{\orbit}[1]{\orbitname_{#1}} 
\newcommand{\hassename}{\mathcal H}     
\newcommand{\glbidx}[2]{
    \hassename_{#2}(#1)}                
\newcommand{\isoidx}[3]{
    \hassename_{#3}(#1|#2)}             
\newcommand{\orbidx}[3]{
    \hassename_{#3}^*(#1|#2)}           

\title[Subalgebras and Hasse Constants]
{Investigation of Subalgebra Lattices by Means of Hasse Constants}
\author{Petr Vojt\v echovsk\'y}
\thanks{\\ While working on this paper the author has been partially
supported by Grant Agency of Charles University, grant number
269/2001/B-MAT/MFF}
\address{Department of Mathematics, Iowa State University, Ames, IA 50011, U.S.A.}
\email{petr@math.du.edu}



\begin{document}


\begin{abstract} Hasse constants and their basic properties are introduced to
facilitate the connection between the lattice of subalgebras of an algebra $C$
and the natural action of the automorphism group $\aut{C}$ on $C$. These
constants are then used to describe the lattice of subloops of the smallest
nonassociative simple Moufang loop.
\end{abstract}

\maketitle

\section{Introduction}\label{Sc:Introduction}

\noindent To completely describe the lattice of subalgebras $\subalg{C}$ of a
finite algebra $C$ is a difficult task. Moreover, it is not obvious how to
store the information about $\subalg{C}$ efficiently, as the cardinality and
complexity of $\subalg{C}$ is typically much larger than that of $C$.
Fortunately, sometimes there is a procedure that allows us to calculate the
join $A\join B$ and meet $A\meet B$ for every $A$, $B\in\subalg{C}$. For
instance, it is easy to find the join and meet in any boolean algebra $C$,
although $|\subalg{C}|$ grows exponentially in $|C|$. It is this ability to
calculate $\meet$ and $\join$ that is often understood as a complete
description of $\subalg{C}$.

Most of the time we are not so lucky, though, and there is no apparent way to
find joins and meets. The main reason is that $A\join B$ and $A\meet B$ can be
far from both $A$ and $B$ in the lattice $\subalg{C}$. It is therefore more
convenient to have access to a procedure that gives a complete \emph{local}
description of $\subalg{C}$. Assuming that it is possible to find all maximal
subalgebras of $A\le C$ and all subalgebras $B\le C$ in which $A$ is maximal,
the lattice $\subalg{C}$ can be built up inductively. We will refer to all
subalgebras immediately above and immediately below $A$ in $\subalg{C}$ as
\emph{neighbors} of $A$, and we denote the set they form by $\nbd{A}$.

In this context, it is worth paying attention to the automorphism group
$\aut{C}$ and its natural action on $C$, since the neighborhoods of $A$ and $B$
will be ``the same'' for $A$, $B\in\subalg{C}$ belonging to the same orbit of
transitivity of $\aut{C}$. Thus, the lattice $\subalg{C}$ can be fully
described as long as we find
\begin{enumerate}
\item[$(\ell_1)$] one representative $A$ from each orbit of $\aut{C}$,
\item[$(\ell_2)$] the neighborhood $\nbd{A}$ for every representative $A$, and
\item[$(\ell_3)$] an automorphism of $C$ mapping $B$ onto $A$, for every
representative $A$ and every $B$ from the orbit of $A$.
\end{enumerate}
To save space, we can store subalgebras by their generating sets, and
substitute $\nbd{A}$ and the automorphisms required by $(\ell_3)$ with an
efficient algorithm producing those.

The purpose of this paper is twofold. First, to introduce a general
tool---Hasse constants---that is of some help in all three tasks $(\ell_1)$,
$(\ell_2)$, $(\ell_3)$. Secondly, to use Hasse constants to describe the
subloop lattice of the smallest nonassociative simple Moufang loop $\paige{2}$.
The investigation of $\subalg{\paige{2}}$ occupies most of this paper, and is
inevitably of rather detailed nature. We maintain that the power of Hasse
constants is sufficiently demonstrated by the fact that $\subalg{\paige{2}}$
was not known before (see Acknowledgement), especially given the importance of
$\paige{2}$ for the real octonions.

Although considerable invention will be required in each particular case, we
believe that Hasse constants will help to keep track in investigation of any
subalgebra lattice.

A word about the notation: we write $\cyclic{n}$ for the cyclic group of order
$n$, $D_n$ for the dihedral group of order $2n$, and $E_{2^n}$ for the
elementary abelian $2$-group of order $2^n$. A subalgebra generated by the set
$S$ well be denoted by $\spn{S}$.

\section{Hasse constants}

\noindent Let $A$, $B$, $C$ be finite (universal) algebras, $A\le C$. For $X\le
C$, let $\orbit{X}$ denote the orbit of $X$ under the natural action of
$\aut{C}$ on the set of subalgebras of $C$ isomorphic to $X$. We will speak of
the subalgebras of $C$ isomorphic to $X$ as \emph{copies} of $X$ in $C$. Define
\begin{eqnarray*}
    \glbidx{B}{C}&=&|\{B_0\le C;\;B_0\cong B\}|,\\
    \isoidx{A}{B}{C}&=&|\{B_0\le C;\;A\le B_0\cong B\}|.
\end{eqnarray*}
Furthermore, when $B\le C$, let
\begin{displaymath}
    \orbidx{A}{B}{C}=|\{B_0\le C;\; A\le B_0,\,B_0\in\orbit{B}\}|.
\end{displaymath}
In words, $\glbidx{B}{C}$ counts the number of copies of $B$ in $C$,
$\isoidx{A}{B}{C}$ counts the number of copies of $B$ in $C$ containing $A$,
and $\orbidx{A}{B}{C}$ counts the number of copies of $B$ in $C$ containing
$A$ and in the same orbit as $B$.

Yet another description of these constants is perhaps the most appealing. For
$B\le C$, the constant $\glbidx{B}{C}$ counts the number of edges connecting
$C$ to a copy of $B$ in the complete Hasse diagram of $\subalg{C}$. The
remaining two constants can be interpreted in a similar way. We will therefore
refer to them jointly as \emph{Hasse constants}.

Note that $\isoidx{A}{B}{C}=\orbidx{A}{B}{C}$ if $\aut{C}$ acts transitively
on the copies of $B$ in $C$.

\begin{lemma}\label{Lm:Counting}
Let $A$, $B$, $C$ be algebras, $A\le C$.
\begin{enumerate}
\item[\textrm{(i)}] If $B'\cong B$, $C'\cong C$, then
    $\glbidx{B}{C}=\glbidx{B'}{C'}$.
\item[\textrm{(ii)}] If $A'\in\orbit{A}$, $B'\cong B$, then
    $\isoidx{A}{B}{C}=\isoidx{A'}{B'}{C}$.
\item[\textrm{(iii)}] If $A'\in\orbit{A}$, $B\le C$, $B'\in\orbit{B}$, then
    $\orbidx{A}{B}{C}=\orbidx{A'}{B'}{C}$.
\end{enumerate}
\begin{proof}
Part (i) is obvious from the definition of $\glbidx{B}{C}$. The equality
$\isoidx{A}{B}{C}=\isoidx{A}{B'}{C}$ holds if $B\cong B'$. Let
$A'\in\orbit{A}$, and let $f\in\aut{C}$ be an automorphism mapping $A$ to $A'$.
Then $\isoidx{A}{B}{C}=\isoidx{f(A)}{f(B)}{f(C)}=\isoidx{A'}{f(B)}{C}
=\isoidx{A'}{B}{C}$. This proves (ii). Part (iii) is similar (use
$\orbit{B}=\orbit{B'}$).
\end{proof}
\end{lemma}

\begin{figure}[ht]
    \centering
    \input{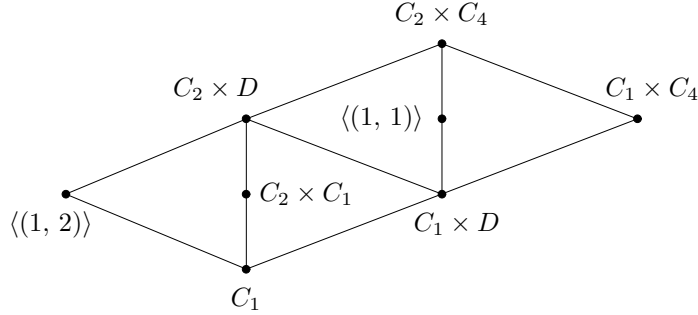}
    \caption[]{Lattice of subgroups of $\cyclic{2}\times\cyclic{4}$}
    \label{Fg:HasseGroup}
\end{figure}

\begin{example}
This example shows that the constants $\isoidx{A}{B}{C}$, $\isoidx{A'}{B}{C}$
may differ even though $A\cong A'$. Let $C$ be the group
$\cyclic{2}\times\cyclic{4}$, $\cyclic{2}=\{0,\,1\}$, $\cyclic{4}=\{0$, $1$,
$2$, $3\}$, and denote by $D=\{0,\,2\}$ the two-element subgroup of
$\cyclic{4}$. The lattice of subgroups of $C$ is depicted in Figure
$\ref{Fg:HasseGroup}$. With $A = \cyclic{1}\times
D\cong\cyclic{2}\times\cyclic{1} = A'$, we have $\isoidx{A}{\cyclic{4}}{C} =
2\ne 0 = \isoidx{A'}{\cyclic{4}}{C}$.
\end{example}

\setlength{\unitlength}{1.00mm}

\begin{proposition}\label{Pr:Counting}
Let $C$ be an algebra, $A$, $B\le C$. Let $A_1$, $\dots$, $A_m$ be
representatives from all orbits $\orbit{A_1}$, $\dots$, $\orbit{A_m}$ of the
action of $\aut{C}$ on the copies of $A$ in $C$. Similarly, let $B_1$, $\dots$,
$B_n$ be representatives for $B$. Then
\begin{eqnarray}
    \isoidx{A}{B}{C}&=&\sum_{j=1}^n{\orbidx{A}{B_j}{C}},\label{Eq:SumOrbits}\\
    \glbidx{A}{B}\cdot|\orbit{B}|&=&\sum_{i=1}^m
        |\orbit{A_i}|\cdot\orbidx{A_i}{B}{C},\label{Eq:StructuralOrb}\\
    \glbidx{A}{B}\cdot\glbidx{B}{C}&=&\sum_{i=1}^m
        |\orbit{A_i}|\cdot\isoidx{A_i}{B}{C}.\label{Eq:StructuralIso}
\end{eqnarray}
When $\aut{C}$ acts transitively on the copies of $B$ $($i.e., when $n=1)$,
then $(\ref{Eq:StructuralOrb})$ coincides with $(\ref{Eq:StructuralIso})$. When
$\aut{C}$ acts transitively on the copies of $A$ $($i.e., when $m=1)$, then
\begin{eqnarray}
    \glbidx{A}{B}\cdot|\orbit{B}|
        =\glbidx{A}{C}\cdot\orbidx{A}{B}{C},\label{Eq:TransitiveOrb}\\
    \glbidx{A}{B}\cdot\glbidx{B}{C}
        =\glbidx{A}{C}\cdot\isoidx{A}{B}{C}.\label{Eq:TransitiveIso}
\end{eqnarray}
\begin{proof}
Since every copy of $B$ in $C$ belongs to exactly one orbit $\orbit{B_j}$,
$(\ref{Eq:SumOrbits})$ follows. To establish $(\ref{Eq:StructuralOrb})$, count
twice the cardinality $t$ of
\begin{displaymath}
    \{(A_0,\,B_0);\;A\cong A_0\le B_0\in\orbit{B}\}.
\end{displaymath}
On the one hand,
\begin{displaymath}
    t=\sum_{B_0\in\orbit{B}}{\glbidx{A}{B_0}}\seq{\ref{Lm:Counting}\mathrm{(i)}}
     \sum_{B_0\in\orbit{B}}{\glbidx{A}{B}}=\glbidx{A}{B}\cdot|\orbit{B}|.
\end{displaymath}
On the other hand,
\begin{displaymath}
    t=\sum_{A_0\le C,\, A_0\cong A}{\orbidx{A_0}{B}{C}}
     = \sum_{i=1}^m\sum_{A_0\in\orbit{A_i}}\orbidx{A_0}{B}{C}
        \seq{\ref{Lm:Counting}\mathrm{(iii)}}
        \sum_{i=1}^m|\orbit{A_i}|\cdot\orbidx{A_i}{B}{C}.
\end{displaymath}
The proof of $(\ref{Eq:StructuralIso})$ is similar to
$(\ref{Eq:StructuralOrb})$. Just count twice the cardinality of the set
\begin{displaymath}
    \{(A_0,\,B_0);\;A\cong A_0\le B_0\le C,\, B_0\cong B\}.
\end{displaymath}
When $m=1$, $(\ref{Eq:TransitiveOrb})$ and $(\ref{Eq:TransitiveIso})$ follow
immediately from $(\ref{Eq:StructuralOrb})$ and $(\ref{Eq:StructuralIso})$,
respectively.
\end{proof}
\end{proposition}

\section{Finite simple Moufang loops and loops of type $\chein$}
\label{Sc:Loops}

\noindent Loops satisfying one of the equivalent \emph{Moufang identities}, for
instance the identity
\begin{equation}\label{Eq:MI}
    ((xy)x)z=x(y(xz)),
\end{equation}
are called \emph{Moufang loops} \cite{Pflugfelder}. By a result of Kunen
\cite{Kunen}, every quasigroup satisfying $(\ref{Eq:MI})$ is a Moufang loop.
Obviously, every group is a Moufang loop. Moufang loops are \emph{power
associative} (i.e., every $1$-generated subloop is a group), in fact
\emph{diassociative} (i.e., every $2$-generated subloop is a group). Every
element $x$ of a Moufang loop has a (unique) two sided inverse $x^{-1}$.

Paige \cite{Paige}, Doro \cite{Doro} and Liebeck \cite{Liebeck} showed that
there is only one class of nonassociative finite simple Moufang loops,
consisting of loops $\paige{q}$, one for each finite field $\field{q}$.

These loops are best studied via composition algebras. Following
\cite{SpringerVeldkamp}, let $\octo{q}$ be the unique split octonion algebra
over $\field{q}$. Then $\paige{q}$ is isomorphic to the multiplicative loop of
elements of norm $1$ in $\octo{q}$ modulo the center. The algebra $\octo{q}$
was first constructed by Zorn as follows. Given a prime power $q$, let $\cdot$
be the standard dot product and $\times$ the standard vector product on
$\field{q}^3$. Then the algebra of \emph{vector matrices}
\begin{displaymath}
    x=\vm{a}{\alpha}{\beta}{b}\;\;\;\;(a,\,b\in \field{q},\,\alpha,\,\beta\in
    \field{q}^3)
\end{displaymath}
with addition defined entry-wise and multiplication governed by
\begin{equation}\label{Eq:Zorn}
    \vm{a}{\alpha}{\beta}{b}\vm{c}{\gamma}{\delta}{d} =
    \vm{ac+\dpr{\alpha}{\delta}}{a\gamma+\alpha d-\vpr{\beta}{\delta}}
            {\beta c+b\delta+\vpr{\alpha}{\gamma}}{\dpr{\beta}{\gamma}+bd}
\end{equation}
is isomorphic to $\octo{q}$. The norm on $\octo{q}$ coincides with the
determinant $\det{x}=ab-\dpr{\alpha}{\beta}$.  The neutral element is
\begin{displaymath}
    \neutral=\vm{1}{(0,\,0,\,0)}{(0,\,0,\,0)}{1},
\end{displaymath}
and every element $x$ with nonzero norm has inverse
\begin{displaymath}
    x^{-1}=\frac{1}{\det{x}}\cdot\vm{b}{-\alpha}{-\beta}{a}.
\end{displaymath}

The order of $\paige{q}$ is $q^3(q^4-1)$ when $q$ is even, and $q^3(q^4-1)/2$
when $q$ is odd \cite{Paige}.

Notably, the loop $\paige{2}$ has also connections to the standard real
(division) octonion algebra. Namely, it is isomorphic to the integral real
octonions of norm $1$ modulo the center (cf. \cite{Coxeter68},
\cite{VojtechovskyJA}).

Let us recall loops of type $\chein=\schein$, first constructed in
\cite{Chein74}. For a group $G$ of order $n$, define new multiplication $\cdot$
on $G\times\cyclic{2}$ by
\begin{displaymath}
    (g,\,i)\cdot(h,\,j)=((g^{(-1)^j}h^{(-1)^{i+j}})^{(-1)^j},\,i+j).
\end{displaymath}
Then $(G\times \cyclic{2},\,\cdot)$ is a Moufang loop, and we denote it by
$\schein$. It is nonassociative if and only if $G$ is nonabelian (cf.\
\cite{Chein74}).

Write $\schein=G\cup Gu$ for some element $u\in\schein\setminus G$. Lemma
\ref{Lm:StructuralM2n} (cf.\ \cite[Prop 4.12]{VojtechovskyThesis}) is easy to
prove once you realize that
\begin{enumerate}
\item[-] every element of $Gu$ is of order $2$,
\item[-] $G\cdot G = Gu\cdot Gu = G$, $G\cdot Gu = Gu\cdot G = Gu$,
\item[-] every subloop $H\not\le G$ of $\schein$ satisfies $|H\cap G|=|H\cap Gu|$.
\end{enumerate}

\begin{lemma}\label{Lm:StructuralM2n}
Let $G$ be a group of order $n$, and let $\schein=G\cup Gu$ be constructed as
above.
\begin{enumerate}
\item[\textrm{(i)}] We have
    \begin{displaymath}
        \glbidx{\cyclic{m}}{\schein}=\left\{\begin{array}{ll}
            \glbidx{\cyclic{m}}{G},&\textrm{if $m\ne 2$,}\\
            \glbidx{\cyclic{m}}{G}+n,&\textrm{if $m=2$.}
        \end{array}\right.
    \end{displaymath}
\item[\textrm{(ii)}] $\spn{H,\,gu}\cong E_{2^{k+1}}$ for every $g\in G$,
    $H\le G$, $H\cong E_{2^k}$, $k\ge 0$.
\item[\textrm{(iii)}] For $k\ge 1$,
    \begin{displaymath}
        \glbidx{E_{2^k}}{\schein}=\left\{\begin{array}{ll}
            0,&\textrm{if $2^{k-1}\nmid n$,}\\
            \glbidx{E_{2^k}}{G}+\glbidx{E_{2^{k-1}}}{G}\cdot
            n\cdot 2^{1-k},&\textrm{otherwise.}
        \end{array}\right.
    \end{displaymath}
\item[\textrm{(iv)}] $\spn{g,\,hu}\cong S_3$ for every $g$, $h\in G$ with
    $\order{g}=3$.
\item[\textrm{(v)}] When $\glbidx{\cyclic{3}}{G}\ne 0$ and $\glbidx{S_3}{G}=0$,
    then $\glbidx{G}{\schein}=1$.
\end{enumerate}
\end{lemma}

It was proved in \cite{VojtechovskyRM} that $\schein$ is presented (in the
variety of Moufang loops) by
\begin{equation}\label{Eq:Pres}
    \spn{x,\,y,\,u;\;R,\,u^2=(xu)^2=(yu)^2=((xy)u)^2=\neutral},
\end{equation}
whenever $G$ is a $2$-generated group with presentation
\begin{displaymath}
    \spn{x,\,y;\;R}.
\end{displaymath}
In particular, presentations for the loops $M(S_3)$ and $M(A_4)$ can be
obtained from this result.

\section{Main goal}

\noindent We proceed to describe the subloop lattice of $\paige{2}$, guided by
steps $(\ell_1)$, $(\ell_2)$, $(\ell_3)$ of Section \ref{Sc:Introduction}. We
fulfil $(\ell_1)$ and give reasonable amount of details with respect to
$(\ell_2)$ and $(\ell_3)$. In particular, we calculate all Hasse constants
$\glbidx{B}{C}$, $\isoidx{A}{B}{C}$, $\orbidx{A}{B}{C}$, for $A<B\le
C=\paige{2}$.

At several places, the reader will be kindly asked to verify a few details by
straightforward, easy calculations. Most of these calculations are reduced to a
quick glance into Table \ref{Tb:Main}. The table itself can be checked for
accuracy within minutes, using Lemma \ref{Lm:Involutions}. No machine
computation is needed.

\section{Possible subloops}\label{Sc:Possible}

\noindent Fix $F=\field{2}$ and $C=\paige{2}$. It is easy to see that $C$
consists of $120$ elements of order $1$, $2$, $3$. More precisely,
\begin{displaymath}
    x=\vm{a}{\alpha}{\beta}{b}
\end{displaymath}
satisfies $\order{x}=2$ if and only if $a=b$ and $x\ne\neutral$; and
$\order{x}=3$ if and only if $a\ne b$. To linearize our notation, we write
$x=\iinv{\alpha}{\beta}{a}$ when $\order{x}=2$, and $x=\tri{\alpha}{\beta}{a}$
when $\order{x}=3$. Since $a$ can be calculated from $\alpha$, $\beta$ when
$\order{x}=2$, we further simplify involutions to $x=\inv{\alpha}{\beta}$. Note
that $\tri{\alpha}{\beta}{a}^{-1}=\tri{\alpha}{\beta}{1+a}$, where the addition
is modulo $2$.

Elementary counting reveals that there are $63$ involutions and $56$ elements
of order $3$ in $C$. Using the language of Hasse constants,
$\glbidx{\cyclic{2}}{C}=63$, $\glbidx{\cyclic{3}}{C}=56/2=28$.

Chein classified all nonassociative Moufang loops of order at most $63$
\cite{Chein78}, and we will call such Moufang loops \emph{small}. Since $C$ has
$120$ elements, every proper subloop of $C$ is small and can be found in
Chein's list.

As in \cite{Pflugfelder}, we say that a finite loop $L$ has the \emph{weak
Cauchy property} when it contains a subloop of order $p$ for every prime $p$
dividing $\order{L}$. It has the \emph{weak Lagrange property} if $|H|$ divides
$|L|$ for every $H\le L$. Finally, $L$ has the \emph{strong Cauchy}
(\emph{Lagrange}) \emph{property} if every subloop of $L$ has the weak Cauchy
(Lagrange) property.

Since $5$ divides $|C|$, $C$ does not have the weak Cauchy property. However,
it follows from \cite[Ch.\ XIV]{Chein78} that all small Moufang loops have it.
(They also have the strong Lagrange property. If one proves that every
$\paige{q}$ has the strong Lagrange property, it will follow that all Moufang
loops have it (cf.\ \cite{CKRV}). As of now, this is an open question.
Glauberman proved \cite{Glauberman} that all Moufang loops of odd order have
the strong Lagrange property.)

\begin{corollary}\label{Cr:RestrictedOrders}
The order of every proper subloop of $C$ is $2^r3^s$, for some $r$, $s$.
\end{corollary}

\begin{lemma}\label{Lm:HasInvolution}
Let $x$, $y\in C$, $\order{x}=\order{y}=3$, $y\not\in\spn{x}$. Then
$\spn{x,\,y}$ contains an involution.
\begin{proof}
We may assume that $x=\tri{\alpha}{\beta}{1}$, $y=\tri{\gamma}{\delta}{1}$, for
some $\alpha$, $\beta$, $\gamma$, $\delta\in F^3$. Then exactly one of the two
elements $xy$, $x^2y$ is of order $2$.
\end{proof}
\end{lemma}

This means that $9$ does not divide the order of any subgroup of $C$. Every
group of order $24$ contains an element of order at least $4$ (the only two
nonabelian groups of order $24$ with Sylow $2$-subgroups isomorphic to $E_8$
are $D_6\times\cyclic{2}$ and $A_4\times\cyclic{2}$). Hence $|G|\in\{1$, $2$,
$3$, $4$, $6$, $8$, $12$, $16$, $32$, $48\}$ for every subgroup $G$ of $C$. It
is not obvious, at least to the author, that $C$ contains no subgroups of order
$16$, necessarily isomorphic to $E_{16}$. It is true, however, and we prove it
in Section \ref{Sc:Lattice}. Hence $|G|=\{1$, $2$, $3$, $4$, $6$, $8$, $12\}$.

Chein concludes in \cite[Ch.\ XII]{Chein78} that every small Moufang loop
containing no element of order greater that $3$ is necessarily of the form
$\schein$ for some nonabelian group $G$. Thanks to the restrictions on $|G|$,
there are only two candidates for $G$, namely $S_3$ (the symmetric group of
order $6$) and $A_4$ (the alternating group of order $12$).

\begin{corollary}\label{Cr:PossibleSubloops}
A nontrivial subloop of $C$ is isomorphic to
\begin{equation}\label{Eq:PossibleSubloops}
    \cyclic{2},\,\cyclic{3},\,E_4,\,S_3,\,E_8,\,A_4,\,
    M(S_3)\text{\ or }M(A_4).
\end{equation}
In particular, $C$ has the strong Lagrange property.
\end{corollary}

All loops listed in $(\ref{Eq:PossibleSubloops})$ actually occur as subloops of
$C$, as we shall see.

\section{Automorphisms}\label{Sc:Autos}

\noindent We construct three kinds of automorphisms of $C$.

\begin{lemma}\label{Lm:LieAuto}
Let $f:F^3\to F^3$ be a nonsingular linear transformation. Define
$\widehat{f}:\octo{2}\to\octo{2}$ by
\begin{displaymath}
    \widehat{f}\vm{a}{\alpha}{\beta}{b}=\vm{a}{f(\alpha)}{f(\beta)}{b}.
\end{displaymath}
Then $\widehat{f}\in\aut{\octo{2}}$ if $($and only if$)$ $f$ is an automorphism
of the Lie algebra $(F^{3},\,+,\,\times)$.
\begin{proof}
Linearity is obvious and the rest follows by straightforward computation using
Zorn's multiplication $(\ref{Eq:Zorn})$.
\end{proof}
\end{lemma}

Identify $\pi\in S_3$ with the linear transformation $F^3\to F^3$, $(\alpha_1$,
$\alpha_2$, $\alpha_3)\mapsto(\alpha_{\pi(1)}$, $\alpha_{\pi(2)}$,
$\alpha_{\pi(3)})$. By Lemma \ref{Lm:LieAuto}, $\widehat{\pi}\in\aut{C}$. To
keep the terminology simple, we will call such automorphisms permutations (of
coordinates).

Define $\dswitch:\octo{2}\to\octo{2}$ by
\begin{equation}\label{Eq:DSwitch}
    \dswitch\vm{a}{\alpha}{\beta}{b}=\vm{b}{\beta}{\alpha}{a},
\end{equation}
and verify that $\dswitch\in\aut{\octo{2}}$.

Finally, we focus on conjugations. Not every conjugation in a Moufang loop is
an automorphism.

\begin{lemma}\label{Lm:Conj}
Let $Q$ be a simple Moufang loop. For $x\in Q$ define $\conj{x}:Q\to Q$ by
$\conj{x}(y)=x^{-1}yx$. Then $\conj{x}$ is a nontrivial automorphism of $Q$ if
and only if $\order{x}=3$.
\begin{proof}
By \cite[Thm IV.1.6]{Pflugfelder}, $\conj{x}$ is a pseudo-automorphism with
companion $x^{-3}$. So $\conj{x}$ is an automorphism whenever $\order{x}$
divides $3$.

Conversely, if $\conj{x}$ is a nontrivial automorphism of $Q$, then it must be
a pseudo-automorphism with companions $\neutral$ and $x^{-3}$. By \cite[Thm
IV.1.8]{Pflugfelder}, the set of all companions of $\conj{x}$ equals
$\neutral\nucleus{Q}$, where $\nucleus{Q}$ is the nucleus of $Q$. Since $Q$ is
simple, we must have $x^{-3}=\neutral$.
\end{proof}
\end{lemma}

\begin{remark}
The conclusion of Lemma \ref{Lm:LieAuto} remains valid over any finite field
$\field{q}$, but we then get $\widehat{-\pi}\in\aut{\octo{q}}$, rather than
$\widehat{\pi}\in\aut{\octo{q}}$. The map $\dswitch:\octo{q}\to\octo{q}$
defined by $(\ref{Eq:DSwitch})$ is an automorphism if and only if $q$ is even.
\end{remark}

\section{Subloops isomorphic to $\cyclic{2}$}

\noindent The detailed discussion concerning $\subalg{C}$ starts here.

\begin{lemma}\label{Lm:Involutions}
Let $x=\iinv{\alpha}{\beta}{n}$, $y=\iinv{\gamma}{\delta}{m}$ be two
involutions, $x\ne y$, and let $z=\tri{\varepsilon}{\varphi}{t}$ be an element
of order $3$ in $C$. Then:
\begin{enumerate}
\item[\textrm{(i)}]
    $\order{xy}=2$ if and only if
    $\spn{x,\,y}\cong E_4$ if and only if
    $\dpr{\alpha}{\delta}=\dpr{\beta}{\gamma}$,
\item[\textrm{(ii)}]
    $\order{xy}=3$ if and only if
    $\spn{x,\,y}\cong S_3$ if and only if
    $\dpr{\alpha}{\delta}\ne\dpr{\beta}{\gamma}$,
\item[\textrm{(iii)}]
    $x$ is contained in a copy of $S_3$,
\item[\textrm{(iv)}]
    every copy of $S_3$ contains an involution of the form
    $\iinv{\fullblank}{\fullblank}{0}$,
\item[\textrm{(v)}]
    $\order{zx}=2$ if and only if
    $\dpr{\alpha}{\varphi}+\dpr{\beta}{\varepsilon}=n$,
\item[\textrm{(vi)}]
    $z$ is contained in a copy of $S_3$.
\end{enumerate}
\begin{proof}
The involution $x$ commutes with $y$ if and only if $\order{xy}=2$. Since
\begin{displaymath}
    xy=\vm{nm+\dpr{\alpha}{\delta}}{\fullblank}{\fullblank}
        {nm+\dpr{\beta}{\gamma}},
\end{displaymath}
parts (i) and (ii) follow. Part (v) is proved similarly.

Let $x=\iinv{\alpha}{\beta}{n}$. Without loss of generality, assume that
$\beta\ne 0$. Pick $\gamma$, $\delta$ so that $\dpr{\alpha}{\delta}=0$,
$\dpr{\beta}{\gamma}\ne 0$. Then choose $m\in\{0,\,1\}$ so that
$y=\iinv{\gamma}{\delta}{m}\in C$. Then $\spn{x,\,y}\cong S_3$, and (iii) is
proved.

Let $G\le C$, $G\cong S_3$, and suppose that $x=\iinv{\alpha}{\beta}{1}$,
$y=\iinv{\gamma}{\delta}{1}\in G$, $x\ne y$. Then
\begin{displaymath}
    xy=\vm{1+\dpr{\alpha}{\delta}}{\alpha+\gamma+\vpr{\beta}{\delta}}
        {\beta+\delta+\vpr{\alpha}{\gamma}}{1+\dpr{\beta}{\gamma}}.
\end{displaymath}
Since $\order{xy}=3$, we have $\dpr{\alpha}{\delta}\ne\dpr{\beta}{\gamma}$. In
other words, $\dpr{\alpha}{\delta}+\dpr{\beta}{\gamma}=1$. Then the third
involution $xyx\in G$ equals
\begin{displaymath}
    \vm{1+\dpr{\alpha}{\delta}+\dpr{(\alpha+\gamma)}{\beta}}
        {\fullblank}{\fullblank}{\fullblank}
    =\vm{\dpr{\alpha}{\beta}}{\fullblank}{\fullblank}{\fullblank}.
\end{displaymath}
Now, $\dpr{\alpha}{\beta}=0$ since $\det{x}=1$, and we are through with (iv).

Since $\det{z}\ne 0$, we can assume that the first coordinate of both
$\varepsilon$ and $\varphi$ is equal to $1$. Then $\alpha=(1,0,0)$,
$\beta=(0,0,0)$ and $n=1$ make $x=\iinv{\alpha}{\beta}{n}$ into an involution
satisfying $\order{zx}=2$, by (v). This proves (vi).
\end{proof}
\end{lemma}

Let us write $\alpha_1\alpha_2\alpha_3$ for the vector $\alpha=(\alpha_1$,
$\alpha_2$, $\alpha_3)\in F^3$, and $\weight{\alpha}$ for
$\alpha_1+\alpha_2+\alpha_3$ (the \emph{weight} of $\alpha$).

Introduce $x_0=\inv{111}{111}$ as the canonical involution of $C$.

Observe that when $x$, $y$ are two involutions of $C$ generating a subgroup
isomorphic to $S_3$ then $\conj{yx}\in\aut{C}$ maps $x$ to $y$.

\begin{proposition}\label{Pr:TransitiveOnC2}
$\aut{C}$ acts transitively on the $63$ copies of $\cyclic{2}$ in $C$.
\begin{proof}
We show how to map an arbitrary $x=\iinv{\alpha}{\beta}{n}$ onto $x_0$. By
Lemma \ref{Lm:Involutions}(iii), $x$ is contained in a copy of $S_3$. Then, by
Lemma \ref{Lm:Involutions}(iv) and the observation immediately preceding this
Proposition, we can assume that $n=0$.

Let $r=\weight{\alpha}$, $s=\weight{\beta}$. Using $\dswitch$ from Section
\ref{Sc:Autos}, we can assume that $r\ge s$. We now fix $y=\inv{100}{100}$ and
proceed to transform $x$ into $x'$ so that $x'=x_0$, or $x'=y$, or
$\spn{x',\,x_0}\cong S_3$, or $\spn{x',\,y}\cong S_3$.

When $r\not\equiv s\pmod{2}$, then $\spn{x,\,x_0}\cong S_3$, by Lemma
\ref{Lm:Involutions}(ii). So assume that $r\equiv s$. Since $n=0$, we have
$s>0$, and thus $(r,\,s)=(1,\,1)$, $(2,\,2)$, $(3,\,1)$, or $(3,\,3)$. Every
permutation of coordinates can be made into an automorphism of $C$, as we have
seen in Section \ref{Sc:Autos}. Moreover, $x_0$ is invariant under all such
permutations. When $(r,\,s)=(1,\,1)$, transform $x$ into $y$. When
$(r,\,s)=(2,\,2)$, transform $x$ into $x'=\inv{110}{011}$, and note that
$\spn{x',\,y}\cong S_3$. When $(r,\,s)=(3,\,1)$, transform $x$ into
$x'=\inv{111}{001}$, and note again that $\spn{x',\,y}\cong S_3$. Finally, when
$(r,\,s)=(3,\,3)$, we have $x=x_0$.

Now, when $\spn{x',\,x_0}\cong S_3$ or $\spn{x',\,y}\cong S_3$, we can permute
the involutions of $C$ so that $x'$ is mapped to $x_0$ or $y$, respectively.
Since $x_0=\conj{\tri{001}{101}{1}}(y)$, we are done.
\end{proof}
\end{proposition}

Note that, in spirit of $(\ell_3)$, the proof of Proposition
\ref{Pr:TransitiveOnC2} tells us how to construct automorphisms mapping
involutions of $C$ onto the representative $x_0$.

\section{Subloops Isomorphic to $\cyclic{3}$ or $S_3$}\label{Sc:C3}

\noindent In this section, we apply Proposition \ref{Pr:Counting} for the first
time. We have taken advantage of the fact that all permutations and $\dswitch$
leave $x_0$ invariant. To proceed further, we need additional automorphisms
with this property.

Consider $v_0=\tri{010}{110}{0}$, $v_1=\tri{001}{101}{0}$, and define
$\ifix:C\to C$ by $\ifix = \conj{v_1^{-1}}\circ\conj{v_0}$. Then
$\ifix\in\aut{C}$, by Lemma \ref{Lm:Conj}, and $\ifix(x_0)=x_0$.

Set $x_1=\inv{110}{100}$, and let $y_0=x_0x_1=\tri{011}{110}{1}$ be the
canonical element of order $3$.

\begin{proposition}\label{Pr:TransitiveOnC3}
$\aut{C}$ acts transitively on the copies of $S_3$ and $\cyclic{3}$.
\begin{proof}
Since $\glbidx{\cyclic{3}}{S_3}=1$ and $\isoidx{\cyclic{3}}{S_3}{C}>0$, by
Lemma \ref{Lm:Involutions}(vi), it suffices to prove that $\aut{C}$ acts
transitively on the copies of $S_3$. Let $G\cong S_3$, $G=\spn{x,\,y}$,
$\order{x}=\order{y}=2$. By Proposition \ref{Pr:TransitiveOnC2}, we can assume
that $x=x_0$. Write $y=\iinv{\alpha}{\beta}{n}$, $r=\weight{\alpha}$,
$s=\weight{\beta}$. By Lemma \ref{Lm:Involutions}, we have $r\not\equiv s$.
Using $\dswitch$, we can assume that $r>s$. We are going to transform $y$ into
$x_1$.

Assume that $n=1$. Then $\dpr{\alpha}{\beta}=0$. Taking permutations of
coordinates and the possible values of $(r,\,s)$ into account, we transform $y$
into one of $x_2=\inv{010}{000}$, $x_3=\inv{011}{100}$, $x_4=\inv{111}{000}$,
$x_5=\inv{111}{101}$. With $\ifix$ as above, check that all of $\ifix(x_2)$,
$\ifix^{-1}(x_3)$, $\ifix(x_4)$ and $\ifix^{-1}(x_5)$ have zeros on the
diagonal.

We may hence assume that $n=0$. Then $(r,\,s)=(2,\,1)$, and we are done by
permuting coordinates.
\end{proof}
\end{proposition}

\begin{lemma}
$\isoidx{\cyclic{2}}{S_3}{C}=16$, $\glbidx{S_3}{C}=336$,
$\isoidx{\cyclic{3}}{S_3}{C}=12$.
\begin{proof}
Pick an involution $x$. By Proposition \ref{Pr:TransitiveOnC3}, the number of
involutions $y$ satisfying $\order{xy}=3$ is independent of $x$. One can then
immediately see with $x=\inv{100}{100}$, say, that there are $32$ such
involutions. As $\glbidx{\cyclic{2}}{S_3}=3$, we get
$\isoidx{\cyclic{2}}{S_3}{C}=16$. Then, by $(\ref{Eq:TransitiveIso})$,
$\glbidx{S_3}{C} = \glbidx{\cyclic{2}}{C} \cdot \isoidx{\cyclic{2}}{S_3}{C}
\cdot \glbidx{\cyclic{2}}{S_3}^{-1} = 336$. Again by
$(\ref{Eq:TransitiveIso})$, $\isoidx{\cyclic{3}}{S_3}{C}
=\glbidx{\cyclic{3}}{S_3} \cdot \glbidx{S_3}{C} \cdot
\glbidx{\cyclic{3}}{C}^{-1} = 12$.
\end{proof}
\end{lemma}

Note that Lemma \ref{Lm:Involutions} allows us to construct all copies of $S_3$
containing $x_0$, and also all copies of $S_3$ containing $y_0$. Note further
that we did not have to resort to local analysis to find the value of
$\isoidx{\cyclic{3}}{S_3}{C}$.

From this moment on, we will pay less attention to $(\ell_2)$ and $(\ell_3)$.

\section{Subloops isomorphic to $A_4$}

\noindent Fix $z_0=\tri{110}{100}{0}$, and recall that
$\glbidx{\cyclic{2}}{A_4}=3$, $\glbidx{\cyclic{3}}{A_4}=4$.

\begin{proposition}\label{Pr:TransitiveOnA4}
$\aut{C}$ acts transitively on the $63$ copies of $A_4$, and\linebreak
$\isoidx{\cyclic{2}}{A_4}{C}=3$.
\begin{proof}
Working in $C$, we have $\spn{G,\,x}\cong S_3$ or $A_4$ for every copy $G$ of
$\cyclic{3}$ and every involution $x$. Since $\isoidx{\cyclic{3}}{S_3}{C}=12$
and $\glbidx{\cyclic{2}}{S_3}=3$, there are $36$ involutions $x$ in $G$ such
that $\spn{G,\,x}\cong S_3$. Thus $\isoidx{\cyclic{3}}{A_4}{C} = (63-36) \cdot
\glbidx{\cyclic{2}}{A_4}^{-1} = 9$. By $(\ref{Eq:TransitiveIso})$,
$\glbidx{A_4}{C} = \glbidx{\cyclic{3}}{C} \cdot \isoidx{\cyclic{3}}{A_4}{C}
\cdot \glbidx{\cyclic{3}}{A_4}^{-1} = 63$.

As for the transitivity, pick $G\cong A_4$. By Proposition
\ref{Pr:TransitiveOnC2}, we can assume that $G=\spn{x_0,\,z}$, for some
$z=\tri{\varepsilon}{\varphi}{t}$ with $r=\weight{\varepsilon}$,
$s=\weight{\varphi}$. Since $\order{x_0z}=3$, we have $r\not\equiv s$, by Lemma
\ref{Lm:Involutions}(v), and may thus assume that $r>s\ge 1$. Then
$(r,\,s)=(2,\,1)$ is the only possibility, and $z$ can be transformed to $z_0$
or $z_0^{-1}$.
\end{proof}
\end{proposition}

Perhaps it would be more natural to look at the copies of $E_4$ now, however,
the Klein subgroups of $C$ are exceptional in the sense that $\aut{C}$ does not
act transitively on them (cf.\ Section \ref{Sc:E4}). We therefore proceed
towards $3$-generated subloops instead.

\section{Subloops isomorphic to $M(S_3)$}

\noindent Have a look at Table \ref{Tb:Main}. It lists all involutions of $C$
and their relation to a few chosen elements of $C$.

\def\cell#1#2#3#4#5#6{\begin{array}{l}#1\ #2\ #3\\#4\ #5\ #6\end{array}}
\begin{table}[ht]
\begin{small}
\begin{displaymath}
\begin{array}{r||l|l|l|l|l|l|l|l}
    {\alpha\setminus\beta}&000&001&010&011&100&101&110&111\\
    \hline\hline
     000    &                           &\cell{3}{ }{ }{3}{ }{ }
            &\cell{3}{ }{ }{2}{ }{ }    &\cell{2}{3}{ }{3}{3}{2}
            &\cell{3}{ }{ }{3}{ }{ }    &\cell{2}{3}{ }{2}{3}{3}
            &\cell{2}{2}{2}{3}{2}{3}    &\cell{3}{ }{ }{2}{ }{ }\\
    \hline
    001     &\cell{3}{ }{ }{3}{ }{ }    &\cell{2}{2}{3}{2}{1}{2}
            &\cell{2}{3}{ }{3}{3}{2}    &\cell{3}{ }{ }{2}{ }{ }
            &\cell{2}{3}{ }{2}{3}{3}    &\cell{3}{ }{ }{3}{ }{ }
            &\cell{3}{ }{ }{2}{ }{ }    &\cell{2}{2}{2}{3}{2}{3}\\
    \hline
    010     &\cell{3}{ }{ }{2}{ }{ }    &\cell{2}{2}{3}{3}{3}{3}
            &\cell{2}{3}{ }{2}{2}{3}    &\cell{3}{ }{ }{3}{ }{ }
            &\cell{2}{3}{ }{3}{2}{2}    &\cell{3}{ }{ }{2}{ }{ }
            &\cell{3}{ }{ }{3}{ }{ }    &\cell{2}{2}{3}{2}{3}{2}\\
    \hline
    011     &\cell{2}{2}{2}{3}{3}{3}    &\cell{3}{ }{ }{2}{ }{ }
            &\cell{3}{ }{ }{3}{ }{ }    &\cell{2}{3}{ }{2}{2}{3}
            &\cell{3}{ }{ }{2}{ }{ }    &\cell{2}{3}{ }{3}{2}{2}
            &\cell{2}{2}{2}{2}{3}{2}    &\cell{3}{ }{ }{3}{ }{ }\\
    \hline
    100     &\cell{3}{ }{ }{3}{ }{ }    &\cell{2}{3}{ }{2}{3}{2}
            &\cell{2}{2}{3}{3}{2}{1}    &\cell{3}{ }{ }{2}{ }{ }
            &\cell{2}{2}{3}{2}{2}{3}    &\cell{3}{ }{ }{3}{ }{ }
            &\cell{3}{ }{ }{2}{ }{ }    &\cell{2}{3}{ }{3}{3}{3}\\
    \hline
    101     &\cell{2}{3}{ }{2}{3}{2}    &\cell{3}{ }{ }{3}{ }{ }
            &\cell{3}{ }{ }{2}{ }{ }    &\cell{2}{2}{3}{3}{2}{2}
            &\cell{3}{ }{ }{3}{ }{ }    &\cell{2}{2}{3}{1}{2}{3}
            &\cell{2}{3}{ }{3}{3}{3}    &\cell{3}{ }{ }{2}{ }{ }\\
    \hline
    110     &\cell{2}{3}{ }{3}{2}{3}    &\cell{3}{ }{ }{2}{ }{ }
            &\cell{3}{ }{ }{3}{ }{ }    &\cell{2}{2}{2}{2}{3}{3}
            &\cell{3}{ }{ }{2}{ }{ }    &\cell{2}{2}{3}{3}{3}{2}
            &\cell{2}{3}{ }{2}{2}{2}    &\cell{3}{ }{ }{3}{ }{ }\\
    \hline
    111     &\cell{3}{ }{ }{2}{ }{ }    &\cell{2}{3}{ }{3}{2}{3}
            &\cell{2}{2}{3}{2}{3}{3}    &\cell{3}{ }{ }{3}{ }{ }
            &\cell{2}{2}{2}{3}{3}{2}    &\cell{3}{ }{ }{2}{ }{ }
            &\cell{3}{ }{ }{3}{ }{ }    &\cell{ }{ }{ }{ }{ }{ }
\end{array}
\end{displaymath}
\end{small}
\caption{Involutions in $\paige{2}$ and their relation to a few elements of
$\paige{2}$. The cell in row $\alpha$ and column $\beta$ corresponds to
involution $x=\inv{\alpha}{\beta}$. It contains the values $\order{x_0x}$,
$\order{x_1x}$, $\order{y_0x}$ (in the first row), $\order{(z_0^{-1}x_0z_0)x}$,
$\order{u_1x}$, $\order{u_2x}$ (in the second row), where $x_0$, $y_0$, $z_0$,
$u_1$ and $u_2$ are as in Theorem \ref{Th:Main}, and $x_1=\inv{110}{100}$. The
values of $\order{x_1x}$, $\order{y_0x}$, $\order{u_1x}$, $\order{u_2x}$ are
calculated only when $\order{x_0x}=2$. Moreover, $\order{y_0x}$ is calculated
only when $\order{x_1x}=2$. No orders are calculated for $e$ and $x_0$.}
\label{Tb:Main}
\end{table}

By \cite{Chein74} or \cite{VojtechovskyRM}, $\glbidx{S_3}{M(S_3)}=3$. By Lemma
\ref{Lm:StructuralM2n}(i), $\glbidx{\cyclic{2}}{M(S_3)}=9$.

Introduce $u_0=\inv{000}{110}$.

\begin{proposition}\label{Pr:TransitiveOnM12}
$\isoidx{S_3}{M(S_3)}{C}=1$. In particular, $\aut{C}$ acts transitively on the
$112$ copies of $M(S_3)$.
\begin{proof}
In view of Proposition \ref{Pr:TransitiveOnC3}, it suffices to prove
$h=\isoidx{S_3}{M(S_3)}{C}=1$ and count the copies of $M(S_3)$.

By $(\ref{Eq:Pres})$, $M(S_3)$ is presented by $\spn{x,\,y,\,u;\;
x^2=y^2=(xy)^3=u^2=(xu)^2=(yu)^2=((xy)u)^2=\neutral}$. Using Lemma
\ref{Lm:Involutions}, verify that $x=x_0$, $y=x_1$ and $u=u_0$ satisfy these
presenting relations, i.e., that $h\ge 1$.

Let $G$ be a copy of $S_3$. By Proposition \ref{Pr:TransitiveOnC3}, we can
assume that $G=\spn{x_0,\,x_1}$. According to Table \ref{Tb:Main}, there are
$6$ involutions $u$ such that $\order{x_0u}=\order{x_1u}=\order{(x_0x_1)u}=2$
(recall that $x_0x_1=y_0$). Since $M_{12}(G)\setminus G$ consists solely of
involutions, we have proved $h\le 1$.

By $(\ref{Eq:TransitiveIso})$, $\glbidx{M(S_3)}{C} = \glbidx{S_3}{C} \cdot
\isoidx{S_3}{M(S_3)}{C} \cdot \glbidx{S_3}{M(S_3)}^{-1} = 112$.
\end{proof}
\end{proposition}

\section{Subloops isomorphic to $M(A_4)$}

\noindent We have $\glbidx{A_4}{M(A_4)}=1$ as a special case of Lemma
\ref{Lm:StructuralM2n}(v).

Introduce $u_1=\inv{001}{001}$.

\begin{proposition}\label{Pr:TransitiveOnMA4}
$\isoidx{A_4}{M(A_4)}{C}=1$. In particular, $\aut{C}$ acts transitively on the
$63$ copies of $M(A_4)$.
\begin{proof}
In view of Proposition \ref{Pr:TransitiveOnA4}, it suffices to prove
$h=\isoidx{A_4}{M(A_4)}{C}=1$ and count the copies of $M(A_4)$.

By $(\ref{Eq:Pres})$, $M(A_4)$ is presented by $\spn{x,\,y,\,u;\; x^2 = y^3 =
(xy)^3 = u^2 = (xu)^2 = (yu)^2 = ((xy)u)^2 = \neutral}$. Verify that $x=x_0$,
$y=z_0$ and $u=u_1$ do the job, hence $h\ge 1$.

Let $G$ be a copy of $A_4$. By Proposition \ref{Pr:TransitiveOnA4}, we can
assume that $G=\spn{x_0,\,z_0}$. Then $v=z_0^{-1}x_0z_0\in G$ is an involution.
According to Table \ref{Tb:Main}, there are $13$ involutions $u$ such that
$\order{x_0u}=\order{vu}=2$. (One of them is $x_0v$.) That is why $h\le 1$.

By $(\ref{Eq:TransitiveIso})$, $\glbidx{M(A_4)}{C} = \glbidx{A_4}{C} \cdot
\isoidx{A_4}{M(A_4)}{C} \cdot \glbidx{A_4}{M(A_4)}^{-1} = 63$.
\end{proof}
\end{proposition}

Let us calculate a few more Hasse constants.

\begin{lemma}\label{Lm:MoreHasse}
We have
\begin{displaymath}
    \begin{array}{l}
    \isoidx{\cyclic{2}}{A_4}{C}=3,\,
    \glbidx{\cyclic{3}}{M(A_4)}=4,\,
    \isoidx{\cyclic{3}}{M(A_4)}{C}=9,\\
    \glbidx{\cyclic{2}}{M(S_3)}=9,\,
    \isoidx{\cyclic{2}}{M(S_3)}{C}=16,\
    \glbidx{\cyclic{2}}{M(A_4)}=15,\\
    \isoidx{\cyclic{2}}{M(A_4)}{C}=15,\,
    \glbidx{\cyclic{3}}{M(S_3)}=1,\
    \isoidx{\cyclic{3}}{M(S_3)}{C}=4,\\
    \glbidx{S_3}{M(A_4)}=16,\,
    \isoidx{S_3}{M(A_4)}{C}=3.
    \end{array}
\end{displaymath}
\begin{proof}
Since $\glbidx{\cyclic{2}}{A_4}=3$, $(\ref{Eq:TransitiveIso})$ yields
$\isoidx{\cyclic{2}}{A_4}{C}=3$. As $\glbidx{\cyclic{m}}{\schein}$ is known
(Lemma \ref{Lm:StructuralM2n}), the value of $\isoidx{\cyclic{m}}{\schein}{C}$
can be calculated by $(\ref{Eq:TransitiveIso})$, too.

It remains to find $\glbidx{S_3}{M(A_4)}$ and $\isoidx{S_3}{M(A_4)}{C}$. Let
$M=M(A_4)=G\cup Gu$, where $G\cong A_4$. Every subgroup of $M$ isomorphic to
$S_3$ can be written as $\spn{x,\,yu}$ for some $x$, $y\in G$, $\order{x}=3$,
and there are exactly six choices of $(x,\,y)$. Since
$\glbidx{\cyclic{3}}{A_4}=4$, we have $\glbidx{S_3}{M(A_4)}=2\cdot 4\cdot
12/6=16$. Consequently, $\isoidx{S_3}{M(A_4)}{C}=3$.
\end{proof}
\end{lemma}

\section{Subloops isomorphic to $E_4$}\label{Sc:E4}

\noindent As announced before, we show that $\aut{C}$ does not act transitively
on the copies of $E_4$.

Introduce $u_2=\inv{100}{010}$.

\begin{lemma}\label{Lm:AtMostTwo}
Let $\spn{x,\,y}$ be one of the $315$ copies of $E_4$ in $C$. Then there is
$\varphi\in\aut{C}$ such that $\varphi(x)=x_0$ and
$\varphi(y)\in\{u_1,\,u_2\}$.
\begin{proof}
Recall that $\isoidx{\cyclic{2}}{S_3}{C}=16$. Therefore, given any involution
$x$, there are $63-1-2\cdot 16=30$ involutions $y$ such that $\spn{x,\,y}\cong
E_4$. Hence, $\glbidx{E_4}{C}=63\cdot 30/(2\cdot 3)=315$.

As always, we may assume that $x=x_0$, $y=\iinv{\alpha}{\beta}{n}$,
$\weight{\alpha}=r$, $\weight{\beta}=s$, $r\equiv s$, and $r\le s$. When
$(r,\,s)=(0,\,2)$, transform $y$ into $u_0$; if $(r,\,s)=(1,\,1)$, into $u_1$
or $u_2$, depending on $n$; if $(r,\,s)=(1,\,3)$, into $u_3=\inv{001}{111}$; if
$(r,\,s)=(2,\,2)$, into $u_4=\inv{110}{110}$ or $u_5=\inv{011}{101}$.

Recall the automorphism $\ifix$ from Section \ref{Sc:C3}, and check that
$\ifix(u_4)=u_1$, $\ifix(u_3)=u_2$, $\ifix(u_5)=u_3$,
$\ifix^{-1}(u_5)=\dswitch(u_0)$. Thus $u_4$ can be transformed into $u_1$, and
each of $u_0$, $u_3$, $u_5$ can be transformed into $u_2$.
\end{proof}
\end{lemma}

Assume, for a while, that $\aut{C}$ acts transitively on the $315$ copies of
$E_4$. Then, by $(\ref{Eq:TransitiveIso})$, $\isoidx{E_4}{A_4}{C} =
\glbidx{E_4}{A_4} \cdot \glbidx{A_4}{C} \cdot \glbidx{E_4}{C}^{-1} = 1\cdot
63/315$, a contradiction. Hence, by Lemma \ref{Lm:AtMostTwo}, there are $2$
orbits of transitivity $\orbitname^+$, $\orbitname^-$, with representatives
$E_4^+=\spn{x_0,\,u_1}$, $E_4^-=\spn{x_0,\,u_2}$.

Since $\glbidx{\cyclic{2}}{E_4}=3$, we have $\orbidx{\cyclic{2}}{E_4^+}{C} =
6/2 = 3$. Then $\orbidx{\cyclic{2}}{E_4^-}{C}=12$. By
$(\ref{Eq:TransitiveOrb})$, $|\orbitname^+| = \glbidx{\cyclic{2}}{C} \cdot
\orbidx{\cyclic{2}}{E_4^+}{C} \cdot \glbidx{\cyclic{2}}{E_4}^{-1} = 63$ and,
similarly, $|\orbitname^-| = 252$.

By $(\ref{Eq:StructuralOrb})$, $63 = \glbidx{E_4}{A_4} \cdot \glbidx{A_4}{C} =
|\orbitname^+| \cdot \isoidx{E_4^+}{A_4}{C} + |\orbitname^-| \cdot
\isoidx{E_4^-}{A_4}{C} = 63 \cdot \isoidx{E_4^+}{A_4}{C} + 252 \cdot
\isoidx{E_4^-}{A_4}{C}$. This is only possible when $\isoidx{E_4^+}{A_4}{C}=1$
and $\isoidx{E_4^-}{A_4}{C}=0$. In other words, a copy of $E_4$ is contained in
$A_4$ if and only if it belongs to $\orbitname^+$.

Let us have a look at the relation between $E_4$ and $M(S_3)$.

\begin{lemma}\label{Lm:UseCosetDistribution}
$\isoidx{E_4^+}{M(S_3)}{C}=0$, $\isoidx{E_4^-}{M(S_3)}{C}=4$.
\begin{proof}
Consider $E_4^+=\spn{x_0,\,u_1}$. Assume that there is $G\cong S_3$ such that
$E_4^+\le \schein$. Since $\{\neutral$, $g_0$, $g_1$, $g_2\}=E_4^+\not\le G$,
there is exactly one involution $g_i$ in $G$, say $g_0$. Write
$g_i=\inv{\alpha_i}{\alpha_i}$ for appropriate vectors $\alpha_i\in F^3$, and
note that $\alpha_0+\alpha_1+\alpha_2=0$.

There is $y=\inv{\gamma}{\delta}\in G$ such that $\spn{y,\,g_0}=G$. Then
$\order{yg_0}=3$, $\order{yg_1}=\order{yg_2}=2$. By Lemma \ref{Lm:Involutions},
$\dpr{\gamma}{\alpha_i}\ne\dpr{\delta}{\alpha_i}$ if and only if $i=0$. Hence
$0 = \dpr{\gamma}(\alpha_0+\alpha_1+\alpha_2) \ne
\dpr{\delta}(\alpha_0+\alpha_1+\alpha_2) = 0$, a contradiction.

Inevitably, $\isoidx{E_4^+}{M(S_3)}{C}=0$. We proceed to calculate
$\isoidx{E_4^-}{M(S_3)}{C}$. Since $\glbidx{E_4}{M(S_3)}=9$, by Lemma
\ref{Lm:StructuralM2n}, we have  $9 \cdot 112 = \glbidx{E_4}{M(S_3)} \cdot
\glbidx{M(S_3)}{C} = |\orbitname^+| \cdot \isoidx{E_4^+}{M(S_3)}{C} +
|\orbitname^-| \cdot \isoidx{E_4^-}{M(S_3)}{C} = 63 \cdot 0+252 \cdot
\isoidx{E_4^-}{M(S_3)}{C}$.
\end{proof}
\end{lemma}

Finally, we have a look at the constants
$c^\varepsilon=\isoidx{E_4^\varepsilon}{M(A_4)}{C}$, for
$\varepsilon\in\{+,\,-\}$.

\begin{lemma}
With the above notation for $c^+$, $c^-$, we have
\begin{enumerate}
\item[\textrm{(i)}] $(c^+,\,c^-)\in\{(3,\,4)$, $(7,\,3)$, $(11,\,2)$, $(15,\,1)$,
    $(19,\,0)\}$,
\item[\textrm{(ii)}] $c^+\le 7$,
\item[\textrm{(iii)}] $c^-\le 3$.
\end{enumerate}
Hence $c^+=7$ and $c^-=3$.
\begin{proof}
Since $\glbidx{\cyclic{2}}{A_4}=3$ and $\glbidx{E_4}{A_4}=1$, we have
$\glbidx{E_4}{M(A_4)}=19$, by Lemma \ref{Lm:StructuralM2n}. Formula
$(\ref{Eq:StructuralIso})$ then yields $19 \cdot 63 = \glbidx{E_4}{M(A_4)}
\cdot \glbidx{M(A_4)}{C} = |\orbitname^+| \cdot c^+ + |\orbitname^-| \cdot c^-
= 63c^+ + 252c^- = (c^+ + 4c^-)\cdot 63$. In particular, $c^++4c^-=19$, and (i)
follows.

Let $E_4^+=\spn{x_0,\,u_1}$. We are trying to find a group $G\cong A_4$ such
that $E_4^+\le\schein$. We look again at the distribution of the $3$
involutions $x_0$, $u_1$, $x_0u_1$ in the cosets $G$, $Gu$. There are two
possibilities: either $E_4^+\le G$, or $|E_4^+\cap G|=2$.

Suppose that $E_4^+\le G$. As $\glbidx{E_4}{A_4}=1$ and
$\glbidx{A_4}{M(A_4)}=1$, there is at most one subloop $M\cong M(A_4)$ such
that $E_4^+\le M$ in such a case.

Now suppose that $|E_4^+\cap G|=2$. Then $E_4^+\cap G$ is one of the three
$2$-element subgroups of $E_4^+$. Let us call it $H$. Since
$\glbidx{\cyclic{2}}{A_4}=3$ and $\glbidx{A_4}{M(A_4)}=1$, there are at most
$3$ subloops $M\cong M_{24}(G)$ such that $H\le G\le M$. Because there are
three ways to choose $H$ in $E_4^+$, there are at most $3\cdot 3=9$ subloops
$M\cong M(A_4)$ such that $E_4^+\le M$.

Altogether, $c^+\le 1+9=10$. By (i), $c^+\le 7$, and (ii) is finished,

Let $E_4^-=\spn{x_0,\,u_2}$. We are trying to find a group $G\cong A_4$ such
that $E_4^-\le M_{24}(G)$. Since $\isoidx{E_4^-}{M(A_4)}{C}=0$, the group
$E_4^-$ is not contained in $G$, i.e., $|E_4^-\cap G|=2$. By Proposition
\ref{Pr:TransitiveOnC2}, we can assume that $E_4^-\cap G=\{\neutral,\,x_0\}$.
If there is such a group $G$, there is also an element
$y=\tri{(\gamma_1,\,\gamma_2,\,\gamma_3)}{(\delta_1,\,\delta_2,\,\delta_3)}{n}$
such that $\spn{x_0,\,y}=G\cong A_4$, i.e.,
\begin{equation}\label{Eq:Aux5}
    \order{yx_0}=3,\,\order{yu_2}=2,\,\order{y(x_0u_2)}=2.
\end{equation}
By Lemma \ref{Lm:Involutions}, the system of equations $(\ref{Eq:Aux5})$ is
equivalent to
\begin{equation}\label{Eq:Aux6}
    \begin{array}{rcl}
        \delta_1+\delta_2+\delta_3+\gamma_1+\gamma_2+\gamma_3&=&1,\\
        \delta_1+\gamma_2&=&1,\\
        \delta_2+\gamma_1&=&1.
    \end{array}
\end{equation}
In particular, $\gamma_3+\delta_3=1$. There are $4$ solutions to
$(\ref{Eq:Aux6})$, namely
\begin{displaymath}
    \binom{\gamma_1,\,\gamma_2,\,\gamma_3}{\delta_1,\,\delta_2,\,\delta_3}
    =\binom{k,\,k+1,\,m}{k,\,k+1,\,m+1},\;\;\;\ (k,\,m=0,\,1).
\end{displaymath}
This is easy to see since both $(\gamma_1,\,\gamma_2)=(0,\,0)$, $(1,\,1)$ lead
to $\det{y}=0$. Hence, there are at most $8$ candidates for $y$ (with $n=0$,
$1$). However, if $\spn{x_0,\,y}$ is isomorphic to $A_4$, then every element of
order $3$ in $\spn{x_0,\,y}$ must satisfy $(\ref{Eq:Aux6})$. There are $8$
elements of order $3$ in $A_4$, and thus there is at most $1$ subloop $M(G)$
satisfying all of our restrictions.

Because our choice of $x_0\in E_4^-\cap G$ was one of three possible choices,
we conclude that $c^-\le 3$.

Combine (i), (ii), (iii) to get $c^+=7$, $c^-=3$.
\end{proof}
\end{lemma}

\section{Subloops isomorphic to $E_8$}\label{Sc:E8}

\noindent Recall the representatives $E_4^+=\spn{x_0,\,u_1}$,
$E_4^-=\spn{x_0,\,u_2}$, and observe that the loop $\spn{x_0,\,u_1,\,u_2}$ is a
group isomorphic to $E_8$.

\begin{lemma}\label{Lm:dMinusdPlus}
$\isoidx{E_4^+}{E_8}{C}=3$, $\isoidx{E_4^-}{E_8}{C}=1$.
\begin{proof}
Write $d^\varepsilon=\isoidx{V_4^\varepsilon}{E_8}{C}$. We have seen that both
$d^+$, $d^-$ are positive. Inspection of Table \ref{Tb:Main} reveals that there
are $12$ involutions $y\not\in E_4^+$ such that $\order{x_0y}=\order{u_1y}=2$.
This immediately shows that $d^+\le 3$. In fact, $y = \inv{000}{110}$,
$\inv{010}{010}$, $\inv{010}{100}$ yield $3$ different copies of $E_8$. Thus
$d^+=3$.

Yet another inspection of Table \ref{Tb:Main} shows that there are $12$
involutions $y\not\in E_4^-$ such that $\order{x_0y}=\order{u_2y}=2$. This
means that $d^-\le 3$, but we prove more. The group $E_4^-$ is contained in $4$
copies of $M(S_3)$, by Lemma \ref{Lm:UseCosetDistribution}. Let $M$ be one of
them. We can assume that $M=G\cup Gu$, where $G\cong S_3$, $x_0\in G$, $u_2=u$.
Since $G\cdot Gu=Gu$, every element $y$ of $Gu$ satisfies $\order{x_0y}=2$. No
involution $y$ of $G$ satisfies $\order{x_0y}=2$. Since $Gu\cdot Gu=G$ and
$\glbidx{\cyclic{2}}{G}=3$, there are $3$ involutions $y\in Gu$ such that
$\order{yu}=2$. One of them is $x_0u$. Altogether, $d^-\le (12-(3-1)\cdot 4)/4
= 1$.
\end{proof}
\end{lemma}

\begin{lemma}\label{Lm:InsertedLater}
Every copy of $E_8$ in $C$ contains a subgroup from $\orbitname^-$.
\begin{proof}
Note that the proof of Lemma \ref{Lm:AtMostTwo} implies that
$\spn{x_0,\,y}\in\orbitname^+$ if and only if $y$ is a permutation of $u_1$ or
$u_4$, i.e., $y$ is one of the $6$ diagonal elements in Table \ref{Tb:Main}.
Let us denote this set by $S$.

Let $E$ be a copy of $E_8$ in $C$. Without loss of generality, $x_0\in E$.
Assume that $\spn{x_0,\,y}\in\orbitname^+$ for every $y\in
E\setminus\{e,\,x_0\}$. Then $E=S\cup\{e,\,x_0\}$. We proceed to show that
$x=\inv{001}{001}$, $y=\inv{100}{100}\in S$ satisfy
$\spn{x,\,y}\in\orbitname^-$.

We have carefully chosen the notation so that $y$ is the same as in the proof
of Proposition \ref{Pr:TransitiveOnC2}. According to the last line of that
proof, $x_0=\gamma_{\tri{001}{101}{1}}(y)$. Using the same automorphism again,
we get $\gamma_{\tri{001}{101}{1}}(x)=\inv{001}{000}$. Hence
$\spn{x,\,y}\in\orbitname^-$.
\end{proof}
\end{lemma}

\begin{proposition}\label{Pr:TransitiveE8}
The group $\aut{C}$ acts transitively on the $63$ copies of $E_8$. Also,
$\isoidx{E_8}{M(A_4)}{C}=3$.
\begin{proof}
Let $E$, $E'$ be two subgroups of $C$ isomorphic to $E_8$. Then there are $G$,
$G'\in\orbitname^-$ such that $G\le E$, $G'\le E'$, by Lemma
\ref{Lm:InsertedLater}. Since $G$, $G'$ belong to the same orbit, there is
$\varphi\in\aut{C}$ mapping $G$ onto $G'$. As $\isoidx{E_4^-}{E_8}{C}=1$,
$\varphi$ must map $E$ onto $E'$.

By $(\ref{Eq:StructuralIso})$ and Lemma \ref{Lm:dMinusdPlus}, $7 \cdot
\glbidx{E_8}{C} = \glbidx{E_4}{E_8} \cdot \glbidx{E_8}{C} = |\orbitname^+|
\cdot d^+ + |\orbitname^-| \cdot d^- = 63 \cdot 3 + 1 \cdot 252 = 441$. Hence,
$\glbidx{E_8}{C}=63$. Consequently, $(\ref{Eq:TransitiveIso})$ yields
$\isoidx{E_8}{M(A_4)}{C} = \glbidx{E_8}{M(A_4)} \cdot \glbidx{M(A_4)}{C} \cdot
\glbidx{E_8}{C}^{-1} = 3$, and we are finished.
\end{proof}
\end{proposition}

\section{Subloop lattice}\label{Sc:Lattice}

\noindent It is about time to show that $C$ contains no copies of $E_{16}$.
Assume that $G\cong E_{16}$ is a subgroup of $C$. By Proposition
\ref{Pr:TransitiveE8}, we can assume that $\spn{x_0,\,u_1,\,u_2}\le G$. Then
there must be at least $8$ involutions $y$ outside $\spn{x_0,\,u_1,\,u_2}$ in
$C$ such that $\order{x_0y}=\order{u_1y}=\order{u_2y}=2$. Previous inspection
of Table \ref{Tb:Main} provided none, a contradiction.

Let us summarize the results about $\paige{2}$ obtained in this paper.

\begin{theorem}\label{Th:Main}
The smallest $120$-element nonassociative simple Moufang loop $C$ satisfies the
strong Lagrange property\index{strong Lagrange property} but not the weak
Cauchy property\index{weak Cauchy property}. The following loops $($and no
other$)$ appear as subloops of $C$: $\{e\}$, $\cyclic{2}$, $\cyclic{3}$, $E_4$,
$S_3$, $E_8$, $A_4$, $M(S_3)$, $M(A_4)$, and $C$.

The automorphism group $\aut{C}$ acts transitively on the copies of each of
these subloops, with the exception of $E_4$. There are two orbits of
transitivity $\orbitname^+$, $\orbitname^-$ for $E_4$. With the notational
conventions introduced in Section $\ref{Sc:Possible}$, we have the following
orbit representatives: $\spn{x_0}$ for $\cyclic{2}$, $\spn{y_0}$ for
$\cyclic{3}$, $\spn{x_0,\,u_1}$ for $\orbitname^+$, $\spn{x_0,\,u_2}$ for
$\orbitname^-$, $\spn{x_0,\,y_0}$ for $S_3$, $\spn{x_0,\,u_1,\,u_2}$ for $E_8$,
$\spn{x_0,\,z_0}$ for $A_4$, $\spn{x_0,\,y_0,\,u_0}$ for $M(S_3)$, and
$\spn{x_0,\,z_0,\,u_1}$ for $M(A_4)$, where $x_0=\inv{111}{111}$,
$y_0=\tri{011}{110}{1}$, $z_0=\tri{110}{100}{0}$, $u_0=\inv{000}{110}$,
$u_1=\inv{001}{001}$, and $u_2=\inv{100}{010}$. The subloop structure and Hasse
constants for $C$ are summarized in Figure $\ref{Fg:Lattice}$.
\end{theorem}

\setlength{\unitlength}{0.85mm}
\begin{figure}
    \centering
    \input{hclatt.lp}
    \caption[]
    {The subloop structure and Hasse constants for
    $\paige{2}$. Two nontrivial representatives $A$, $B$ are connected by an
    edge if and only if $\isoidx{A}{B}{\paige{2}}>0$. If $A=\{e\}$ or
    $B=\paige{2}$, the two representatives $A$, $B$ are connected by an edge
    if and only if a copy of $A$ is maximal in $B$. The edge connecting $A$
    and $B$ is thick if and only if a copy of $A$ is maximal in $B$. The
    constants $|\orbit{A}|$, $\glbidx{A}{B}$, $\orbidx{A}{B}{\paige{2}}$
    are located in the diagram as follows: $|\orbit{A}|$ next to $A$;
    $\glbidx{A}{B}$ and $\orbidx{A}{B}{\paige{2}}$ in the box on the edge
    connecting $A$ and $B$, separated by colon.}
    \label{Fg:Lattice}
\end{figure}

\section{Acknowledgement}

\noindent This paper is based on \cite[Ch.\ 5]{VojtechovskyThesis}. Shortly
after I finished working on \cite{VojtechovskyThesis}, Orin Chein brought my
attention to the work of Merlini Guiliani and Polcino Milies
\cite{GuilianiMilies}. In \cite{GuilianiMilies}, the authors studied the
subloop lattice of $C=\paige{2}=GLL(F_2)$ for the first time, and
\begin{enumerate}
\item[-] determined $\glbidx{\cyclic{2}}{C}$, $\glbidx{\cyclic{3}}{C}$, and
    estimated $\isoidx{\cyclic{2}}{E_4}{C}$,
\item[-] used a result analogous to Lemma \ref{Lm:Involutions}(i), (ii),
\item[-] found all possible isomorphism types of subloops of $C$ and listed one
    example for each type, thus establishing the strong Lagrange property for
    $C$,
\item[-] sketched the relations between isomorphism types of subloops of $C$,
\end{enumerate}
all without proofs. They did not notice that $E_8$ is not a maximal subloop of
$C$. I acknowledge that I compared some of my results with
\cite{GuilianiMilies}.

The notion of Hasse constants is new, to my knowledge. The name itself was
suggested to me by Jonathan D.~H.~Smith.

I would also like to thank the reviewer who pointed out that the values of
$\isoidx{\cyclic{3}}{M(A_4)}{C}$ and $\isoidx{\cyclic{2}}{M(S_3)}{C}$ were in
error in the previous version of the paper.

\bibliographystyle{plain}

\end{document}